\documentclass[11pt]{article}

\usepackage{graphicx}
\usepackage{pdflscape}
\usepackage{algorithmic}
\usepackage{algorithm}
\RequirePackage[OT1]{fontenc}
\RequirePackage{amsthm,amsmath}

\numberwithin{equation}{section}
\theoremstyle{plain}
\newtheorem{thm}{Theorem}

\newtheorem{lem}{Lemma}

\newcommand{\bX}{{\bf X}}

\newcommand{\bp}{{\bf p}}

\newcommand{\ba}{{\bf a}}

\newcommand{\br}{{\bf r}}
\newcommand{\bL}{{\bf L}}
\newcommand{\wht}{\widehat}
\newcommand{\wtd}{\widetilde}

\newcommand{\cI}{{\cal I}}
\newcommand{\wbox}{\sqcap\llap{$\sqcup$}}
\newcommand{\bZ}{{\bf Z}}

\newcommand{\btau}{\boldsymbol{\tau}}

\begin{document}

\begin{center}
{\bf MINIMIZING CHANGE-POINT ESTIMATION ERROR}

\medskip
by Chan Hock Peng

{\it National University of Singapore} 

\begin{abstract}
In this paper we consider change-points in multiple sequences with the objective of minimizing the estimation error 
of a sequence by making use of information from other sequences. 
This is in contrast to recent interest on change-points in multiple sequences where the focus is on detection of common
change-points.
We start with the canonical case of a single sequence with constant change-point intensities.
We consider two measures of a change-point algorithm.
The first is the probability of estimating the change-point with no error.
The second is the expected distance between the true and estimated change-points.
We provide a theoretical upper bound for the no error probability, 
and a lower bound for the expected distance, 
that must be satisfied by all algorithms.
We propose a scan-CUSUM algorithm that achieves the no error upper bound and come close to the distance lower bound.
We next consider the case of non-constant intensities and establish sharp conditions under which estimation error can go to zero.
We propose an extension of the scan-CUSUM algorithm for a non-constant intensity function, 
and show that it achieves asymptotically zero error at the boundary of the zero-error regime.
We illustrate an application of the scan-CUSUM algorithm on multiple sequences sharing an unknown, 
non-constant intensity function. 
We estimate the intensity function from the change-point profile likelihoods of all sequences 
and apply scan-CUSUM on the estimated intensity function.
\end{abstract}
\end{center}

\section{Introduction}

This paper considers change-point estimation with a Bayesian framework for the generation of change-points.
Representative papers on the Bayesian model for change-point estimation include Yao (1984), 
Barry and Hartigan (1993), Chib (1998), Lai and Xing (2011), Ko, Chong and Ghosh (2015) and Du, Kao and Kou (2016).
Yao (1984) considered a hidden Markov model (HMM) for change-points with normal observations, 
and showed how the posterior change-point probabilities can be computed using backward-forward formulas.
Lai and Xing (2011) extended the HMM approach 
to multi-parameter exponential families and proposed a bounded complexity mixture to limit computational complexity.
Ko, Chong and Ghosh (2015) estimated the change-point intensities using a Beta prior instead of assuming that they are known.
Du, Kao and Kou (2016) showed that change-points can be estimated consistently using marginal likelihoods of the HMM.

A recent topic of interest in change-point detection is the study of detectability when the data consists of
multiple sequences and change-points are common in a sparse
fraction of the sequences.
Zhang, Siegmund, Ji and Li (2010), motivated by common change-points in multiple DNA sequences of copy number variants,
proposed scan and segmentation algorithms to pool signals using sum of chi-squared statistics.
Mei (2010) proposed a sum of CUSUM test statistic for on-line detection of change-points in multiple data streams and showed that it detects optimally when the fraction of sequences undergo distribution change is correctly specified.
Xie and Siegmund (2013) proposed a generalized likelihood test that is optimal and adaptive to the fraction of sequences undergoing change.
 
Cho and Fryzlewicz (2015) considered the detection of change-points on high-dimensional time-series 
using Wild Binary Segmentation (WBS), 
with thresholding of CUSUM scores from individual sequences.
Jeng, Cai and Li (2010) characterized the regimes under which weak signals from a small fraction of sequences can be pooled for successful detection and showed that higher-criticism test statistics are able to detect at the boundary of such regimes.
Chan and Walther (2015) highlighted how the sequence to signal ratio affect the level of sparsity for which detection is possible.
Optimal detection theories for sparse change-points, including on time-series or with spatial dependence,
were developed in Horv\'{a}th and Huskova (2014), Jirak (2015), Cho (2016), Wang and Samworth (2018), 
Enikeeva and Harchaoui (2019), Pilliart, Carpentier and Verzelen (2020) and Liu, Gao and Samworth~(2021).

\subsection{Contributions of this paper}

Our work involves multiple sequences but our emphasis is different from earlier works 
in that we want to minimize estimation error by applying information from the other sequences.
Key to the method is the proposal of a common change-point intensity function 
and the estimation of this function by using change-point profile likelihoods. 

We provide an upper bound for the no error probability of a change-point algorithm under a constant intensity assumption
and propose a scan-CUSUM algorithm that achieves this bound.
We also provide a lower bound for the expected distance between the true and estimated change-points and show that
our algorithm has an expected distance that is numerically close to this bound.
In the case of a non-constant intensity function we establish sharp conditions 
under which asymptotically zero estimation error is possible, 
and propose an extension of the scan-CUSUM algorithm 
that achieves asymptotically zero error at the boundary of the zero-error regime.
The scan-CUSUM algorithm forms the basis of using information sharing among sequence to improve estimation of change-points
in multiple sequences.

\subsection{Layout of the paper}

In Section 2 we describe the change-point model and provide an upper bound for the probability of estimating the 
change-points with no error and a lower bound for the expected distance between the true and estimated change-points, 
when the change-point intensities are constant. 
In Section 3 we propose a scan-CUSUM algorithm that achieves the no error upper bound and come close 
numerically to the distance lower bound.
The algorithm first localizes a change-point using scan statistics and applies the CUSUM test statistics of WBS
to finalize the change-point estimates.
In Section 4 we extend to non-constant intensities and show that asymptotically zero estimation is not possible outside a zero-error regime.
We extend the scan-CUSUM algorithm to handle change-point intensities that vary with location and show that it achieves 
asymptotically zero error in the zero-error regime.
In Section 5 we illustrate how the scan-CUSUM algorithm can be applied to reduce estimation error 
when there are multiple sequences with change-points generated from a common intensity function.
In Section 6 we perform experiments to compare estimating the change-points one sequence at a time 
versus using all available information to estimate the change-points.
In Section 7 we prove the theoretical results.
In the Appendix we provide additional technical arguments and calculations.

\subsection{Notations}

Let $\#A$ be the number of elements in a set $A$ and let $\emptyset$ be the empty set.
Let $\phi(z) = \tfrac{1}{\sqrt{2 \pi}} e^{-\frac{z^2}{2}}$ be the density, 
and $\Phi(z) = \int_{-\infty}^z \phi(x) dx$ the distribution function, 
of the standard normal.
Let $\lfloor \cdot \rfloor$ be the greatest integer function and $\lceil \cdot \rceil$ the least integer function.
For a number $t$ and set $A$ let $d(t,A) = \min_{u \in A} |u-t|$.
Let ${\bf 1}$ be the indicator function and let ${\bf Z}$ be the set of all integers.
Let $a=o(b)$ if $\tfrac{a}{b} \rightarrow 0$.

\section{Change-point model and estimation error bounds}

Let $\bX = (X(1), \ldots, X(T))$ be the observations of a sequence and let $\mu(t)$ be the mean of $X(t)$.

\medskip
{\sc Change-point model}. 
Assume that there exists random variables $\Delta(t)$ and independent Bernoulli random variables $Y(t)$ such that
\begin{equation} \label{M1}
\mu(t+1) = \mu(t)+ Y(t) \Delta(t). 
\end{equation}
We assume that for all $t$,
$\Delta(t)|Y(t)=1$ follows a common distribution $F$ that has no point mass at zero.
The change-points are $\btau = \{ t: Y(t)=1 \}$.
Let $a(t) = EY(t)$ be the change-point intensity at $t$.
In this section we consider 
\begin{equation} \label{aaq}
a(1)= \cdots = a(T-1)=q \mbox{ for some } q>0.
\end{equation} 

We assume that conditioned on $\mu(t)$, the observations $X(t)$ are independent, 
with
\begin{equation} \label{M2}
X(t) = \mu(t) + \epsilon(t), \quad \epsilon(t) \sim {\rm N}(0,\sigma_X^2).
\end{equation}

\medskip
{\sc Example}.
Yao (1984) considered a HMM with 
\begin{equation} \label{M3}
\mu(t+1) = (1-Y(t)) \mu(t)+ Y(t) \xi(t),
\end{equation}
with $\mu(1), \xi(1), \ldots, \xi(T-1) \sim_{\rm i.i.d.} {\rm N}(0,\sigma_{\xi}^2)$ and independent of all the $Y(t)$.
Conditioned on $Y(t)=1$, $\Delta(t) = \xi(t)-\mu(t)$ is distributed as ${\rm N}(0, 2 \sigma_{\xi}^2)$.
Note however that $\Delta(1), \ldots, \Delta(T-1)$ are dependent.
Yao (1984) showed that under (\ref{M2}) and (\ref{M3}), 
backward-forward formulas can be used to compute $P(Y(t)=1|\bX)$ efficiently.

\medskip
{\sc Remarks}. 
Let $\wht{\btau} = \{ t: \wht Y(t)=1 \}$ be the estimated change-points of a change-point algorithm.
Assuming that the change-point model allows us to compute $P(Y(t)=1|\bX)$ efficiently, 
we maximize $\sum_{t=1}^{T-1} P(\wht Y(t)=Y(t))$ by letting
$\wht \tau = \{ t: P(Y(t)=1|\bX) \geq \tfrac{1}{2} \}$.
However this estimator is overly conservative as we are penalized twice if a change-point is estimated at an incorrect location,
compared to not estimating the change-point at all.
Moreover it does not take into account how close the estimated change-point is to the true change-point.

\medskip
{\sc Notations}. 
\begin{enumerate}
\item Let $J=\# \btau$ be the number of change-points.
Arrange the change-points in $\btau$ as
$$\tau_1 < \cdots < \tau_J.
$$
For completeness define $\tau_0=0$ and $\tau_{J+1}=T$.

\item Let $\wht J = \# \wht{\btau}$ be the number of change-points estimated by an algorithm.
\end{enumerate}

\medskip
In the definitions of $\alpha_T, \beta_T$ and $\gamma_T(\Delta_0)$ below, 
$\tau$ refers to the random selection of one of the $J$ change-points, 
conditioned on $J>0$.

\medskip
{\sc Definitions}.

\begin{enumerate}
\item For each $j$ let 
\begin{equation} \label{kappa}
\kappa(\tau_j) = \# \{ \wht \tau \in \wht{\btau}: |\tau_j-\wht \tau| < \tfrac{1}{2} \min(\tau_{j+1}-\tau_j,\tau_j-\tau_{j-1}) \}.
\end{equation}
\item Let $\alpha_T = P(\kappa(\tau) \neq 1| J \geq 1)$.

\item The probability of estimating a change-point with no error is 
$$\beta_T = P(\wht Y(\tau)=1 \mbox{ and } \kappa(\tau)=1 | J \geq 1).
$$

\item Let $\Delta_0 > 0$.
The expected $L_1$ location error is 
$$\gamma_T(\Delta_0) = E(d(\tau,\wht{\btau})| J \geq 1, \kappa(\tau)=1, |\Delta(\tau)| \geq \Delta_0).
$$
\end{enumerate}

\medskip
{\sc Remarks}. 

\begin{enumerate}
\item In (\ref{kappa}) we assign each $\wht \tau$ to at most one $\tau_j$ 
and $\kappa(\tau_j)$ is the number of $\wht \tau$ that $\tau_j$ is assigned to.
We include $\kappa(\tau)=1$ in the evaluation of a change-point algorithm, 
in $\beta_T$ and $\gamma_T(\Delta_0)$, 
to prevent artificial improvements by either having many estimated change-points 
or avoiding the estimation of change-points with weak signals.

\item We condition on $|\Delta(\tau)| \geq \Delta_0$, 
in the definition of $\gamma_T(\Delta_0)$, 
because $d(\tau,\wht{\btau})$ scales like $[\Delta(\tau)]^{-2}$, 
so the expected $L_1$ location error is infinite without this restriction.
We show in Theorem \ref{thm1} that under (\ref{M1})--(\ref{M2}),
any algorithm subject to the constraint $\alpha_T \rightarrow 0$ as $T \rightarrow \infty$ satisfies 
$$\liminf_{T \rightarrow \infty} \gamma_T(\Delta_0) \geq \gamma_{\rm lower}(\tfrac{\Delta_0}{\sigma_X})
$$
for some positive $\gamma_{\rm lower}(\tfrac{\Delta_0}{\sigma_X})$,
which we characterize by using the following terminologies.
\end{enumerate}

\bigskip \bigskip \bigskip \bigskip \bigskip
{\sc Definitions}.
\begin{enumerate}
\item For a non-negative function $\br = (r(i): i \in {\bf Z})$ such that $s = \sum_{i=-\infty}^{\infty} r(i)$ is positive
and finite, 
${\rm med}(\br)$ is the integer $u$ satisfying
$$\sum_{i=-\infty}^{u-1} r(i) < \tfrac{s}{2} \leq \sum_{i=-\infty}^u r(i),
$$
and ${\rm mode}(\br)$ is the smallest integer $u$ satisfying $r(u) = \max_{i \in \bZ} r(i)$.

\item Let $Z^+(u)$ and $Z^-(u)$ be independent standard normal random variables
and let $\Delta \sim F$ be independent of $Z^+(u)$ and $Z^-(u)$.
Let $S^+(i) = \sum_{u=1}^i Z^+(u)$, $S^-(i) = \sum_{u=1}^i Z^-(u)$ and 
let $\bp = (p(i): i \in {\bf Z})$ 
be such that for $i>0$, 
\begin{eqnarray} \label{pi}
p(i) & = & \exp(\Delta S^+(i)-\tfrac{i \Delta^2}{2}), \\ \nonumber
p(-i) & = & \exp(\Delta S^-(i)-\tfrac{i \Delta^2}{2}), \\ \nonumber
p(0) & = & 1. 
\end{eqnarray}
Define
\begin{equation} \label{glower}
\gamma_{\rm lower}(\Delta_0) = E \Big (|{\rm med}(\bp)| \Big| |\Delta| \geq \Delta_0 \Big).
\end{equation}
It is shown in Lemma \ref{lem1} in Section 7 that $\gamma_{\rm lower}(\Delta_0)<\infty$ for all $\Delta_0>0$.

\item Let
$\Delta \sim F$. Define
\begin{equation} \label{bupper}
\beta_{\rm upper} = \tfrac{1}{2} E(\Delta^2 \nu(\Delta)),
\end{equation}
where $\nu(\Delta) = 2 \Delta^{-2} \exp(-2 \sum_{i=1}^{\infty} i^{-1} \Phi(-\tfrac{\sqrt{i}}{2} |\Delta|))$
is the overshoot function, see (4.37) of Siegmund (1985).
The overshoot function is bounded above by 1 and $\nu(\Delta) \rightarrow 1$ as $\Delta \rightarrow 0$.
\end{enumerate}

\begin{thm} \label{thm1}
Consider $q \rightarrow 0$ as $T \rightarrow \infty$.
Under {\rm (\ref{M1})--(\ref{M2})} all change-point algorithms satisfy
$$\limsup_{T \rightarrow \infty} \beta_T \leq \beta_{\rm upper}.
$$
Let $\Delta_0>0$.
If a change-point algorithm satisfies $\alpha_T \rightarrow 0$ then
$$\liminf_{T \rightarrow \infty} \gamma_T(\Delta_0) \geq \gamma_{\rm lower}(\tfrac{\Delta_0}{\sigma_X}).
$$
\end{thm}

\section{Optimality for constant intensity functions}

The scan-CUSUM algorithm achieves $\beta_{\rm upper}$ when the change-point intensities $a(t)$ are constant.
The algorithm scans windows of observations to localize change-points and applies CUSUM scores within each highlighted window 
to estimate the exact location of change-points.

The use of scan statistics to localize change-points was employed by the Screening and Ranking algorithm of Niu and Zhang (2012).
The use of CUSUM statistics to estimate change-points was used by the WBS algorithm of Fryzlewicz (2014).
The use of scan statistics of multiple windows lengths in scan-CUSUM is motivated by the multiscale methods of
Arias-Castro, Donoho and Huo (2005, 2006).

\medskip
{\sc Definitions}. Let $S(t) = \sum_{u=1}^t X(u)$.
\begin{enumerate}
\item Scan statistics.
For a given window length $2 \ell$ define
$$Z_{\ell}(t) = \sqrt{\tfrac{1}{2 \ell \sigma_X^2}}[S(t+\ell)+S(t-\ell)-2 S(t)],
$$
for $\ell \leq t \leq T-\ell$.

\item CUSUM statistics.
For a given interval $I=(u,v)$ define
\begin{equation} \label{ZI}
Z_I(t) = \sqrt{\tfrac{(v-t)(t-u)}{(v-u) \sigma_X^2}}(\tfrac{S(v)-S(t)}{v-t}-\tfrac{S(t)-S(u)}{t-u}),
\end{equation}
for $u < t < v$.

\item Profile likelihood.
For a given interval $I=(u,v)$ define $\bL(I) = (L_I(t): t \in {\bf Z})$ with
\begin{equation} \label{profile}
L_I(t) = \exp(\tfrac{1}{2} Z_I^2(t)) \mbox{ for } t \in I,
\end{equation}
and $L_I(t)=0$ for $t \not\in I$.
Check that for $t \in I$,
$$L_I(t) = C_I \Big[ \sup_{\mu} \prod_{r=u+1}^t \phi(\tfrac{X_r-\mu}{\sigma_X}) \Big] 
\Big[ \sup_{\nu} \prod_{r=t+1}^v \phi(\tfrac{X_r-\nu}{\sigma_X}) \Big],
$$
where $C_I=(2 \pi)^{\frac{v-u}{2}} \exp[\tfrac{1}{2 \sigma_X^2} \sum_{r=u+1}^v (X_r - \bar X_I)^2]$,
$\bar X_I = \tfrac{1}{v-u} \sum_{r=u+1}^v X_r$.
\end{enumerate}

\medskip
{\sc Remarks}.
The scan-CUSUM algorithm selects $t = \wht \tau$ to maximize $P(\tau=t|\tau \in I)$ once a change-point has been localized
to lie within an interval $I$,
by treating the profile likelihood as the true likelihood.
When the change-point intensity is constant this is equivalent to selecting $t$ to maximize the profile likelihood $L_I(t)$.
However there is no need to refer to the profile likelihood as we are just selecting $t$ to maximize $|Z_I(t)|$.

The profile likelihood function is needed in the following situations.

\begin{enumerate}
\item When the intensity function is constant 
and the objective is to minimize the expected distance between $\tau$ and $\wht \tau$.
Here the optimal estimator selects $\wht \tau$ to be the median of the profile likelihood.
However as the scan-CUSUM has an expected distance close to the universal lower bound, 
we do not explore this further in this paper.

\item When the intensity function is not constant, the maximization of $P(\tau=t|\tau \in I)$ is achieved by selecting
$t=\wht \tau$ to maximize $a(t) L_I(t)$.
We investigate this in Section 4.
\end{enumerate}

\medskip
{\sc Algorithm}. 
Let $1 < \rho \leq 2$ and define $\ell_b = \lceil \rho^b \rceil$.
Write $Z_b(t)$ instead of $Z_{\ell_b}(t)$.
Let $b_T$ be the largest $b$ satisfying $2 \ell_b \leq T-1$.

\begin{algorithm}[h]
\caption{Scan-CUSUM}
\begin{algorithmic}
\STATE Input: $c_{\rm scan} > 0$
\STATE Initialization: $\wht Y(t)=0$ for $1 \leq t \leq T-1$
\STATE For $b=0, \ldots, b_T$:
\STATE \hspace{0.25cm} $A = \{ t: \ell_b \leq t \leq T-\ell_b, \sum_{u=t-\ell_b+1}^{t+\ell_b-1} \wht Y(u)=0 \}$
\STATE \hspace{0.25cm} While $(\max_{t \in A} |Z_b(t)| \geq c_{\rm scan})$
\STATE \hspace{0.5cm} $t^* = {\rm arg \ max}_{t \in A} |Z_b(t)|$
\STATE \hspace{0.5cm} $I = (t^*-\ell_b,t^*+\ell_b)$
\STATE \hspace{0.5cm} $\wht \tau = {\rm arg \ max}_{t \in I} |Z_I(t)|$
\STATE \hspace{0.5cm} Update $\wht Y(\wht \tau) = 1$ 
\STATE \hspace{0.5cm} Update $A$ according to the change in $\wht Y(\wht \tau)$
\STATE Output: $\wht{\btau} = \{ t: \wht Y(t)=1 \}$
\end{algorithmic}
\end{algorithm}

\begin{thm} \label{thm2}
Consider $q=o(\tfrac{1}{\log T})$ as $T \rightarrow \infty$.
Under {\rm (\ref{M1})--(\ref{M2})},
the scan-CUSUM algorithm,
with $c_{\rm scan} = \sqrt{2 \log(T \log T)}$,
satisfies $\alpha_T \rightarrow 0$ and
\begin{equation} \label{thm2eq}
\lim_{T \rightarrow \infty} \beta_T = \beta_{\rm upper}.
\end{equation}
\end{thm}

\bigskip \bigskip \bigskip \bigskip \bigskip \bigskip \bigskip \bigskip \bigskip \bigskip
{\sc Remarks}. 
\begin{enumerate}
\item The threshold $c_{\rm scan} = \sqrt{2 \log (T \log T)}$ is a theoretical value 
chosen to highlight the limiting behavior of the scan-CUSUM.
In practice we may want to apply thresholds to satisfy other criteria,
for example a Type I error probability constraint.

\item Let $\bp$ be as defined in (\ref{pi}) and let $\Delta_0 > 0$.
By the calculations used to show (\ref{thm2eq}),
details in Section 7.2.3,
\begin{equation} \label{gscan}
\lim_{T \rightarrow \infty} \gamma_T(\Delta_0) = \gamma_{\rm scan}(\tfrac{\Delta_0}{\sigma_X}),
\end{equation}
where
$$\gamma_{\rm scan}(\Delta_0) = E \Big( |{\rm mode}(\bp)| \Big| |\Delta| \geq \Delta_0 \Big).
$$
Our numerical investigations below show however that though $\gamma_{\rm scan}(\Delta_0) > \gamma_{\rm lower}(\Delta_0)$,
the difference between $\gamma_{\rm scan}(\Delta_0)$ and $\gamma_{\rm lower}(\Delta_0)$ is small,
so scan-CUSUM works well in terms of minimizing $L_1$ location error.
\end{enumerate}

\medskip
{\sc Numerical comparison}.
Express
\begin{eqnarray*}
\gamma_{\rm lower}(\Delta_0) & = & [1-F(\Delta_0)+F(-\Delta_0)]^{-1} \int_{|\Delta| \geq \Delta_0}
g_{\rm lower}(\Delta) F(d \Delta), \cr
\gamma_{\rm scan}(\Delta_0) & = & [1-F(\Delta_0)+F(-\Delta_0)]^{-1} \int_{|\Delta| \geq \Delta_0}
g_{\rm scan}(\Delta) F(d \Delta), 
\end{eqnarray*}
where
$$g_{\rm lower}(\Delta) = E \Big( |{\rm med}(\bp)| \Big| \Delta \Big) \mbox{ and } 
g_{\rm scan}(\Delta) = E \Big( |{\rm mode}(\bp)| \Big| \Delta \Big).
$$
Table 1,
based on the generation of 10,000 normal sequences,
shows that $g_{\rm scan}(\Delta)$ is about 7--9\% larger than $g_{\rm lower}(\Delta)$.
So the scan-CUSUM algorithm has $L_1$ location error that is close to the lower bound.

\begin{table}[h]
\begin{center}
\begin{tabular}{lcc}
$\Delta$ & $\Delta^2 g_{\rm lower}(\Delta)$ & $\Delta^2 g_{\rm scan}(\Delta)$ \cr \hline
1.0 & 2.71$\pm$0.04 & 2.91$\pm$0.04 \cr
0.3 & 2.75$\pm$0.04 & 2.95$\pm$0.04 \cr
0.1 & 2.78$\pm$0.04 & 2.98$\pm$0.04 \cr
0.03 & 2.82$\pm$0.04 & 3.06$\pm$0.04
\end{tabular}
\end{center}
\caption{ Comparison of $g_{\rm lower}(\Delta)$ and $g_{\rm scan}(\Delta)$.}
\end{table}

\section{Optimality for non-constant intensity functions}

{\sc Change-Point Intensity Model}.
Assume that there exists a distribution $G_q$, 
with mean $q$ and support on $[0,1]$, such that
\begin{equation} \label{at}
a(t) \sim_{\rm i.i.d.} G_q \mbox{ for } 1 \leq t \leq T-1.
\end{equation}
Conditioned on $\{ a(t) \}_{t=1}^{T-1}$, $Y(t)$ are independent Bernoulli($a(t)$) random variables.

\medskip
{\sc Remarks}.
We provide in Theorem \ref{thm3} conditions under which asymptotically zero estimation error is not possible,
assuming that the change-point algorithms are applied with knowledge of the intensity function.
These conditions are sharp because when not satisfied,
we are able to achieve asymptotically zero estimation error.

\begin{thm} \label{thm3}
Assume {\rm (\ref{M1})}, {\rm (\ref{M2})} and {\rm (\ref{at})} with $q \rightarrow 0$ as $T \rightarrow \infty$.
If
\begin{equation} \label{VG}
\limsup_{q \rightarrow 0} \tfrac{1}{q} \int_0^1 [1-G_q(x)]^2 dx > 0,
\end{equation}
then no change-point algorithm is able to achieve $\beta_T \rightarrow 1$.
\end{thm}

{\sc Examples}. 
\begin{enumerate}
\item If $G_q(yq) = G_1(y)$ and $G_1$ is non-negative and bounded with mean 1,
then by a change of variables $y=x/q$,
$$\tfrac{1}{q} \int_0^1 [1-G_q(x)]^2 dx = \int_0^{1/q} [1-G_1(y)]^2 dy
$$
and (\ref{VG}) holds.

\item If $G_q$ is the Beta($\tfrac{q}{1-q},1$) distribution,
then $G_q(x) = x^{\frac{q}{1-q}}$ and as $q \rightarrow 0$,
$$\tfrac{1}{q} \int_0^1 [1-G_q(x)]^2 dx = \tfrac{2q}{1+q} \rightarrow 0,
$$
so (\ref{VG}) does not hold.
\end{enumerate}

\medskip
{\sc Remarks}.
We define below an extended scan-CUSUM algorithm which selects $t=\wht \tau$ to maximize
$$a(t) L_I(t) \Big[ = a(t) \exp(\tfrac{1}{2} Z_I^2(t)) \Big],
$$
where $\bL(I)$ is the profile likelihood (\ref{profile}).
In contrast the scan-CUSUM algorithm for constant intensity function, 
in Section 3, 
selects $t=\wht \tau$ to maximize the CUSUM score $|Z_I(t)|$.
The extended scan-CUSUM achieves $\beta_T \rightarrow 1$ when (\ref{VG}) does not hold.

\medskip
{\sc Algorithm}.
Extended scan-CUSUM.
Proceed as in the scan-CUSUM algorithm.
Replace ``$\wht \tau = {\rm arg  \ max}_{t \in I} |Z_I(t)|$'' in Line 8 of the algorithm by 
``$\wht \tau = {\rm arg  \ max}_{t \in I} a(t) L_I(t)$''.

\medskip
\begin{thm} \label{thm4}
Consider $q=o(\tfrac{1}{\log T})$ as $T \rightarrow \infty$.
Under {\rm (\ref{M1})}, {\rm (\ref{M2})} and {\rm (\ref{at})},
if 
$$\tfrac{1}{q} \int_0^{\infty} [1-G_q(x)]^2 dx \rightarrow 0 \mbox{ as } q \rightarrow 0,
$$
then the extended scan-CUSUM algorithm,
with $c_{\rm scan} = \sqrt{2 \log (T \log T)}$, achieves $\beta_T \rightarrow 1$ as $T \rightarrow \infty$.
\end{thm}

\section{Estimation with multiple sequences}

We show here how scan-CUSUM can be applied to reduce change-point estimation error,
when there are multiple sequences.

\medskip
{\sc Model and notations}.
\begin{enumerate}
\item Let $\bX^n = (X^n(1), \ldots, X^n(T))$ be the $n$th observed sequence.
Assume that there exists unknown intensities $a(t)$ such that
$$Y^n(t) \sim_{\rm indep.} {\rm Bernoulli}(a(t)), 
$$
for $1 \leq t \leq T-1$ and $1 \leq n \leq N$.

\item Let $\btau^n = \{ t: Y^n(t)=1 \}$ be the change-points of the $n$th sequence and let 
$J^n = \# \btau^n$ be the number of change-points.
Arrange the change-points in $\btau^n$ as
$$\tau_1^n < \cdots < \tau_{J^n}^n.
$$

\item Let $\mu^n(t)$ be the mean of $X^n(t)$.
Assume that for each $n$,
$$\mu^n(t+1) = \mu^n(t) + \Delta^n(t) Y^n(t),
$$
with $\Delta^n(t)| Y^n(t)=1$ non-zero.

\item Assume that conditioned on $\mu^n(t)$, 
for $1 \leq t \leq T$ and $1 \leq n \leq N$,
$$X^n(t) = \mu^n(t) + \epsilon^n(t),
$$
with $\epsilon^n(t) \sim_{\rm i.i.d.} {\rm N}(0,\sigma_X^2)$.
\end{enumerate}

\medskip
{\sc Estimation procedure for multiple sequences}.

\begin{enumerate}
\item Let $\wht{\btau}^n$ be the estimated change-points of the $n$th sequence, 
obtained by applying the scan-CUSUM algorithm
in Section 3 on $\bX^n$.

\item Let $\wht J^n = \# \wht{\btau}^n$ be the number of estimated change-points and let $I^n(1), \ldots, I^n(\wht J^n)$
be the intervals identified by the scan-CUSUM algorithm to contain change-points.
Let $\bL^n(I^n(1)), \ldots, \bL^n(I^n(\wht J^n))$ be the corresponding profile likelihoods. 
To simplify notations express the likelihoods as $\bL^n(1), \ldots, \bL^n(\wht J^n)$.

\item Let $\ba = (a(t): 1 \leq t \leq t-1)$ be the true change-point intensity function and let
$\wht a(t)$, $1 \leq t \leq T-1$,
be the MLE of
\begin{equation} \label{ella}
\ell(\ba) = \sum_{n=1}^N \sum_{j=1}^{\hat J^n} \log \Big( \sum_{t \in I^n(j)} a(t) L^n_j(t) \Big) - N \sum_{t=1}^{T-1} a(t).
\end{equation}

\item Let $\wtd{\btau}^n= \{ \wtd \tau_j^n: 1 \leq j \leq \wht J^n \}$ be the change-point estimates,
based on information from all sequences,
with 
$$\wtd \tau_j^n = {\rm arg \ max}_{t \in I^n(j)} \wht a(t) L_j^n(t).
$$
\end{enumerate}

\medskip
The intensity likelihood (\ref{ella}) is derived using an incomplete-data argument, 
with the change-points as latent variables.
This leads to the following EM procedure for computation of $\wht a(t)$.
The technical details are in Appendix~A.

\medskip
{\sc EM procedure to compute $\wht a(t)$}.
We apply the following procedure to estimate $\wht a(t)$ iteratively.

\begin{enumerate}
\item Initialize with $a_0(t) = \tfrac{1}{N(T-1)} \sum_{n=1}^N \wht J^n$ for all $t$.

\item For $k=0,\ldots,K-1$:
Let 
$$a_{k+1}(t) = \tfrac{1}{N} \sum_{(n,j): t \in I^n(j)} \Big( \tfrac{a_k(t) L_j^n(t)}
{\sum_{u \in I^n(j)} a_k(u) L_j^n(u)} \Big).
$$

\item Let $\wht a(t) = a_K(t)$.
\end{enumerate}

\section{Numerical comparisons}

We compare here scan-CUSUM applied in two ways.

\begin{enumerate}
\item No information sharing.
Scan-CUSUM is applied separately on each sequence,
with no estimation of $a(t)$.
The estimates are $\wht{\btau}^n$ in step 1 of the estimation procedure in Section 5.

\item With information sharing.
The intensity function is estimated using all the sequences,
and the estimated intensities are applied in the change-point estimation.
The estimates are $\wtd{\btau}^n$ in step 5 of the estimation procedure in Section 5.
\end{enumerate}

We perform three sets of experiments,
for $N=100$, $T=10$,000 and $q=0.0001$,
each with a different intensity function generator $G_q$.

\begin{enumerate}
\item $G_q$=prob mass 1 at $q$.
Change-points between sequence are independent.
Some loss of accuracy is expected when applying estimated intensities.

\item $G_q$=prob mass 0.01 at $100q$ and 0.99 at 0.
$$( \# \{ n: Y^n(t)=1 \} | a(t)>0 ) \sim_{\rm approx} \mbox{ Poisson}(100Nq=1).
$$
There is some alignment of change-points due to the variability of $a(t)$ 
however it is not strong enough for estimation error to be close to zero.

\item $G_q$=Beta($\tfrac{q}{1-q}$,1).
$$( a(t) | Y^n(t)=1) \sim {\rm Beta}(\tfrac{1}{1-q},1) \doteq {\rm Uniform}(0,1).
$$
The alignment of change-points is strong enough for estimation error to be close to zero.
\end{enumerate}

In each set of experiment we generate $a(t)$, 
$1 \leq t \leq T-1$, 
i.i.d. from $G_q$.
For each $1 \leq n \leq N$, 
we generate $X^n(t)$ using the HMM (\ref{M2}) and (\ref{M3}) with $\sigma_X=\sigma_{\xi}=1$.
We select threshold $c_{\rm scan}=5.05$ for the scan-CUSUM algorithm.
This is based on a Type I error probability of 0.05,
i.e. $P(\wht J^n >0|J^n=0)=0.05$. 

For sequences $n$ with $J^n>0$ we compute
\begin{eqnarray*}
\wht \alpha_T^n & = & \# \{ j: \kappa(\tau_j^n) \neq 1 \}/J^n, \cr
\wht \beta_T^n & = & \# \{ j: \wht Y(\tau_j^n)=1 \mbox{ and } \kappa(\tau_j^n)=1 \}/J^n.
\end{eqnarray*}
We average them over sequence to get $\wht \alpha_T$ and $\wht \beta_T$.

The experiments are repeated 100 times and standard errors are calculated.
Table 2 shows significant improvements when information sharing is applied on aligned change-points. 

\begin{table}[h]
\begin{center}
\begin{tabular}{llll}
$G_q$ & & No info sharing & With info sharing \cr \hline
Prob mass 1 at $q$ & $\alpha_T$ & 0.069$\pm$0.003 & 0.062$\pm$0.002 \cr
& $\beta_T$ & 0.308$\pm$0.006 & 0.271$\pm$0.005 \cr \hline
Prob mass 0.01 at $100q$ & $\alpha_T$ & 0.069$\pm$0.003 & 0.061$\pm$0.002 \cr
and 0.99 at 0 & $\beta_T$ & 0.303$\pm$0.005 & 0.403$\pm$0.006 \cr \hline
Beta($\tfrac{q}{1-q}$,1) & $\alpha_T$ & 0.064$\pm$0.005 & 0.056$\pm$0.005 \cr
& $\beta_T$ & 0.305$\pm$0.010 & 0.835$\pm$0.011 \cr 
\end{tabular}
\end{center}
\caption{ Estimation accuracy comparison for no information sharing versus information sharing.}
\end{table}

\section{Proofs}

By considering $\tfrac{X(t)}{\sigma_X}$ instead of $X(t)$ we may assume without loss of generality that $\sigma_X=1$.
The following lemmas are proved in Appendix B.
Lemma \ref{lem1} ensures $\gamma_{\rm lower}(\Delta_0) < \infty$.
Lemma \ref{lem2} shows that convergence in likelihood functions lead to convergence in their medians.

\medskip
{\sc Definition}.
Given a non-negative function $\br=(r(i): -\infty < i < \infty)$ with $r(0)>0$, 
define
$$W(\br) = \sum_{i=-\infty}^{\infty}|i| \min(1,\tfrac{r(i)}{r(0)}).
$$

\begin{lem} \label{lem1}
{\rm (a)} For any non-negative function $\br$ with $r(0)>0$ and $s=\sum_{i=-\infty}^{\infty} r(i) < \infty$, 
$$|{\rm med}(\br)| \leq 2 W(\br).
$$ 

\medskip \noindent
{\rm (b)} For $\bp$ defined in {\rm (\ref{pi})} and $\Delta_0 > 0$,
$E ( W(\bp) {\bf 1}_{\{ |\Delta| \geq \Delta_0 \}}) < \infty$.
\end{lem}

\begin{lem} \label{lem2}
For non-negative functions $\br$ and $\br_T$ such that $s = \sum_{t=-\infty}^{\infty} r(t) < \infty$
and $\sum_{t=-\infty}^{\infty} |r_T(t)-r(t)| \rightarrow 0$ as $T \rightarrow \infty$,
if
\begin{equation} \label{lem2.1}
\sum_{t=-\infty}^{{\rm med}(\br)-1} r(t) < \tfrac{s}{2} < \sum_{t=-\infty}^{{\rm med}(\br)} r(t),
\end{equation}
 then ${\rm med}(\br_T) = {\rm med}(\br)$ for $T$ large.
\end{lem}

\medskip
{\sc Notations}.
Let $a_T  =O(b_T)$ if $\sup_T |\tfrac{a_T}{b_T}| < \infty$.
Let $X_T = o_p(Y_T)$ if $P(|\tfrac{X_T}{Y_T}| \geq \epsilon) \rightarrow 0$ for all $\epsilon>0$.
Let $X_T=O_p(Y_T)$ if for all $\epsilon > 0$ there exists constants $M_{\epsilon} >0$ 
such that $P(|\tfrac{X_T}{Y_T}| > M_{\epsilon}) \leq \epsilon$.
We apply the convention $\prod_{t=u}^v \cdot=1$ when $u>v$.

\subsection{Proof of Theorem \ref{thm1}}

{\sc Proof of} $\limsup \beta_T \leq \beta_{\rm upper}$.
Let $\tau=\tau_j$ be a randomly selected change-point.
Let $U > 0$ and let $k$ be such that 
$\tau \in I_k: = \{ (k-1)U+1, \ldots, kU \}$.
Let $K=\lfloor \tfrac{T-1}{U} \rfloor$.
Let $J_k$ be the number of change-points in $I_k$.

For $t \in I_k$ let likelihood
\begin{equation} \label{Lkt}
L_k(t) = \prod_{u=(k-1)U+2}^t \phi(X_u-\mu(\tau)) \prod_{u=t+1}^{kU} \phi(X_u-\mu(\tau+1)).
\end{equation}
Consider an oracle estimator
\begin{equation} \label{hatj}
\wht \tau_j = \left\{ \begin{array}{ll} \tau & \mbox{ if } J_{k-1}+J_k+J_{k+1}>1 \mbox{ or } k>K, \cr
{\rm arg \ max}_{t \in I_k} L_k(t) & \mbox{ otherwise.}
\end{array} \right.
\end{equation}
Since $q \rightarrow 0$, the first case in (\ref{hatj}) occurs with prob$\rightarrow 0$.
In the analysis below we condition on the second case.

Let $Z(r) = X(r)-\mu(\tau)$ for $r \leq \tau$, 
$Z(r) = \mu(\tau+1)-X(r)$ for $r > \tau$ and $\Delta = \Delta(\tau) = \mu(\tau+1)-\mu(\tau)$.
Let $S^+(i) = \sum_{u=\tau+1}^{\tau+i} Z(u)$ and $S^-(i) = \sum_{u=\tau-i+1}^{\tau} Z(u)$.

For $t \in I_k$,
\begin{equation} \label{Ltt}
L_k(t) = L_k(\tau) \times \left\{ \begin{array}{ll}
\exp ( \Delta S^+(t-\tau)-\tfrac{(t-\tau) \Delta^2}{2}) & \mbox{ if } t>\tau, \cr
\exp( \Delta S^-(\tau-t)-\tfrac{(\tau-t) \Delta^2}{2}) & \mbox{ if } t < \tau.
\end{array} \right.
\end{equation}
Hence the prob that $\wht \tau_j=\tau$ is 
\begin{eqnarray*}
& & P(\Delta S^+(i) < \tfrac{i \Delta^2}{2} \mbox{ for } 1 \leq i \leq kU-\tau) \cr
& \times & P(\Delta S^-(i) < \tfrac{i \Delta^2}{2} \mbox{ for } 1 \leq i \leq \tau-(k-1)U-1).
\end{eqnarray*}

It is not possible for a change-point estimator to have smaller $\beta_T$ than the oracle estimator hence
\begin{equation} \label{lsb}
\limsup \beta_T \leq P \Big( \max_{1 \leq i \leq V-1} (\Delta S^+(i)-\tfrac{i \Delta^2}{2}) < 0,
\max_{1 \leq i \leq U-V} (\Delta S^-(i)-\tfrac{i \Delta^2}{2}) < 0 \Big),
\end{equation}
where $V$ is uniformly distributed on $\{ 1, \ldots, U \}$ and independent of $S^+(i)$ and $S^-(i)$.

Since $V \uparrow \infty$ and $(U-V) \uparrow \infty$ in prob as $U \rightarrow \infty$,
by (\ref{lsb}),
\begin{eqnarray*}
\limsup \beta_T & \leq & P^2 \Big( \sup_{i \geq 1} (\Delta S^+(i)-\tfrac{i \Delta^2}{2}) < 0 \Big) \cr
& = & E \Big[\exp \Big(-2 \sum_{i=1}^{\infty} i^{-1} \Phi(-\tfrac{\sqrt{i} |\Delta|}{2}) \Big) \Big] = \beta_{\rm upper},
\end{eqnarray*}
see Corollary 8.44 of Siegmund (1985) for the first equality.
$\wbox$

\medskip
{\sc Proof of} $\liminf \gamma_T(\Delta_0) \geq \gamma_{\rm lower}(\Delta_0)$.
Let $\tau=\tau_j$ be a randomly chosen change-point.
Let $\Delta=\Delta(\tau)$.
Since $\alpha_T \rightarrow 0$ it suffices to show that for any $U>0$,
there exists $\gamma_{{\rm lower},U}(\Delta_0)$ such that
\begin{eqnarray} \label{supgamma}
& & \liminf [\gamma_T(\Delta_0) P(|\Delta| \geq \Delta_0)+U \alpha_T] \\ \nonumber
& = & \liminf [E(d(\tau,\wht{\btau}){\bf 1}_{\{ \kappa(\tau)=1 \ {\rm and} \ |\Delta| \geq \Delta_0 \}})
+U \alpha_T] \\ \nonumber 
& \geq & \gamma_{{\rm lower}, U}(\Delta_0) P(|\Delta| \geq \Delta_0),
\end{eqnarray}
and that $\gamma_{{\rm lower},U}(\Delta) \rightarrow \gamma_{\rm lower}(\Delta_0)$ as $U \rightarrow \infty$.

Let $k$ be such that $\tau \in I_k:= \{ (k-1)U+1, \ldots, kU \}$.
Let $K=\lfloor \tfrac{T-1}{U} \rfloor$.
Let $J_k$ be the number of change-points in $I_k$.
Consider an oracle estimator
$$\wht \tau_j = \left\{ \begin{array}{ll} \tau & \mbox{ if } J_{k-1}+J_k+J_{k+1}>1 \mbox{ or } k>K, \cr
{\rm med}(\bL_k) & \mbox{ otherwise,} \end{array} \right.
$$
where $\bL_k= (L_k(t): -\infty < t < \infty)$ with $L_k(t)$ the likelihood (\ref{Lkt}) when $t \in I_k$
and $L_k(t)=0$ otherwise.

By (\ref{Ltt}),
$${\rm med}(\bL_k)-\tau = {\rm med}(\bp_{(k-1)U+1-\tau,kU-\tau}),
$$
where $p_{a,b}(i) = p(i) {\bf 1}_{\{ a \leq i \leq b \}}$ for all $i$. 

It is not possible for a change-point estimator to have smaller
$$E(d(\tau,\wht{\btau}) {\bf 1}_{\{ \kappa(\tau)=1 \mbox{ and } |\Delta| \geq \Delta_0 \}})+U \alpha_T
$$
than the oracle estimator hence the inequality in (\ref{supgamma}) holds with 
\begin{equation} \label{inf2}
\gamma_{{\rm lower},U}(\Delta_0) = E \Big( |{\rm med}(\bp_{V-U,V-1})| \Big| |\Delta| \geq \Delta_0 \Big),
\end{equation}
with $V$ uniformly distributed on $\{ 1,\ldots,U \}$ and independent of $\bp$.

Since $V \uparrow \infty$ and $(V-U) \uparrow \infty$ in prob as $U \rightarrow \infty$,
by Lemmas \ref{lem1}, \ref{lem2} and DOM,
$$\gamma_{{\rm lower},U}(\Delta_0) \rightarrow 
E \Big( |\mbox{med}(\bp)| \Big| |\Delta| \geq \Delta_0 \Big) = 
\gamma_{\rm lower}(\Delta_0)
$$
as $U \rightarrow \infty$. $\wbox$

\subsection{Proof of Theorem \ref{thm2} and (\ref{gscan})}

\subsubsection{Proof of $\alpha_T \rightarrow 0$}

Express $q=\tfrac{\omega^3}{\log T}$,
with $\omega \rightarrow 0$ as $T \rightarrow \infty$.
Select $\tau_j$ randomly from $\btau$.
With prob$\rightarrow$ 1,
\begin{equation} \label{DT}
|\Delta(\tau_j)| \geq \omega, \ |\Delta(\tau_{j-1})| \geq \omega \mbox{ (if } j>1) \mbox{ and }
|\Delta(\tau_{j+1})| \geq \omega \mbox{ (if } j <J).
\end{equation}
Since $q=o(\tfrac{\omega^2}{\log T})$, with prob$\rightarrow 1$,
\begin{eqnarray} \label{taujj}
\tau_{j-1} & \leq & \tau_j - 48 \omega^{-2} \log T, \\ \nonumber
\tau_{j+1} & \geq & \tau_j + 48 \omega^{-2} \log T.
\end{eqnarray}

Under (\ref{DT}),
if $6 \omega^{-2} \log T \leq \ell_b < 12 \omega^{-2} \log T$ then for $i=j-1,j$ and $j+1$,
$Z_b(\tau_i) \sim {\rm N}(\nu_i,1)$ with
$$|\nu_i| = \sqrt{\tfrac{\ell_b}{2}} |\Delta(\tau_j)| \geq \sqrt{3 \log T}.
$$
Hence with prob$\rightarrow 1$,
$|Z_b(\tau_i)| \geq c_{\rm scan}$,
implying that $|\tau_i - \wht \tau_i| \leq 12 \omega^{-2} \log T$.
In view of (\ref{taujj}),
we conclude
\begin{eqnarray} \label{tauh}
\tau_j - \wht \tau_{j-1} & > & \tfrac{1}{2} (\tau_j-\tau_{j-1}), \\ \nonumber
\wht \tau_{j+1} - \tau_j & > & \tfrac{1}{2} (\tau_{j+1}-\tau_j), \\ \nonumber
|\wht \tau_j - \tau_j| & < & \tfrac{1}{2} \min(\tau_{j+1}-\tau_j,\tau_j-\tau_{j-1}).
\end{eqnarray}

With prob$\rightarrow 1$ all $\wht \tau$ in $\btau$ are of the form $\wht \tau_j$, since 
$$\sum_{b=0}^{b_T} \sum_{t=\ell_b}^{T-\ell_b} P ( |Z_b(t)| \geq c_{\rm scan},
Y(u)=0 \mbox{ for } |u-t| < \ell_b) \leq (b_T+1) \Phi(-c_{\rm scan}) \rightarrow 0.
$$
Hence by (\ref{tauh}), 
$P(\kappa(\tau_j)=1) \rightarrow 1$.

\subsubsection{Proof of $\beta_T \rightarrow \beta_{\rm upper}$}

Select $\tau = \tau_j$ randomly from $\btau$ and let $\Delta=\Delta(\tau)$.
Let $\delta>0$ and let
\begin{equation} \label{AT}
A_T = \{ \delta \leq |\Delta| \leq \delta^{-1}, Y(t)=0 \mbox{ for } 1 \leq |t-\tau| \leq 12 \delta^{-2} \log T \}.
\end{equation}
Since $q=o(\tfrac{1}{\log T})$,
\begin{equation} \label{PAT}
P(A_T) \rightarrow P(\delta \leq |\Delta| \leq \delta^{-1}) \mbox{ as } T \rightarrow \infty.
\end{equation}

Let $\cI_T$ be the set of all intervals $I=(t^*-\ell_b,t^*+\ell_b)[=(u^*,v^*)]$ such that
\begin{eqnarray} \label{t/4}
|t^*-\tau| & \leq & (\log T)^{\frac{1}{4}}, \\  \label{show0}
2 \delta^2 \log T < \ell_b & < & 12 \delta^{-2} \log T.
\end{eqnarray}
By Lemma \ref{lem6} in Appendix C,
under $A_T$,
with prob$\rightarrow 1$, the interval $I$ (in Line 7 of the scan-CUSUM algorithm), 
for which $\tau \in I$, 
is in $\cI_T$. 

Let $S^+(i)=\sum_{u=1}^i Z^+(u)$ and $S^-(i)=\sum_{u=1}^i Z^-(u)$,
where $Z^+(u) = \mu(\tau+1)-X(\tau+u)$ and $Z^-(u)=X(\tau-u+1)-\mu(\tau)$.
By Corollary 8.44 of Siegmund (1985),
\begin{equation} \label{Delta2}
\tfrac{1}{2} \Delta^2 \nu(\Delta) = P \Big( \sup_{i \geq 1} (S^+(i)-\tfrac{i|\Delta|}{2})<0,
\sup_{i \geq 1} (S^-(i)-\tfrac{i|\Delta|}{2}) < 0 \Big| \Delta \Big).
\end{equation}

We show below that
\begin{eqnarray} \label{p+}
\sup \Big| \sqrt{2 \ell_b}[Z_I(\tau+i)-Z_I(\tau)]
-(2S^+(i)-i \Delta) \Big| {\bf 1}_{A_T} & \stackrel{p}{\rightarrow} & 0, \\ \label{p-}
\sup \Big| \sqrt{2 \ell_b}[Z_I(\tau-i)-Z_I(\tau)]
-(2S^-(i)-i \Delta) \Big| {\bf 1}_{A_T} & \stackrel{p}{\rightarrow} & 0, 
\end{eqnarray}
with the supremum over $I \in \cI_T$ and $1 \leq i \leq (\log T)^{\frac{1}{4}}$.
We show in Lemma~\ref{lem7} in Appendix C that
\begin{equation} \label{AB}
\sum_{I \in \cI_T} \sum_{u \in I: |u-\tau| \geq (\log T)^{\frac{1}{4}}} P(|Z_I(u)| \geq |Z_I(\tau)|, A_T) \rightarrow 0.
\end{equation}

By (\ref{p+})--(\ref{AB}), the prob that $\wht \tau=\tau$ is asymptotically the prob that
$\sup_{i \geq 1} (S^+(i)-\tfrac{i |\Delta|}{2}) < 0$ and $\sup_{i \geq 1} (S^-(i)-\tfrac{i |\Delta|}{2}) < 0$.
We conclude $\beta_T \rightarrow \beta_{\rm upper}$ from (\ref{Delta2}) and letting $\delta \rightarrow 0$.
$\wbox$

\medskip
{\sc Proofs of} (\ref{p+}) and (\ref{p-}).
Consider $t=\tau+i$ with $1 \leq i \leq (\log T)^{\frac{1}{4}}$.
Since $v^*-u^*=2 \ell_b$, 
we can express
\begin{eqnarray} \label{711a}
\sqrt{2 \ell_b} Z_I(\tau) & = & \lambda(\tau)(\wht \nu-\wht \mu), \\ \nonumber
\mbox{where } \wht \nu & = & \tfrac{S(v^*)-S(\tau)}{v^*-\tau}, \\ \nonumber
\wht \mu & = & \tfrac{S(\tau)-S(u^*)}{\tau-u^*}, \\ \nonumber
\lambda(t) & = & \sqrt{(v^*-t)(t-u^*)}.
\end{eqnarray}
Since $|t^*-\tau| \leq (\log T)^{\frac{1}{4}}$,
\begin{eqnarray} \label{expand1}
\lambda(\tau+i) & = & \sqrt{\ell_b^2-|\tau+i-t^*|^2} \\ \nonumber
& = & \ell_b[1+O((\log T)^{-\frac{3}{2}})], \\ \label{expand2}
\tfrac{v^*-\tau}{v^*-\tau-i} & = & (1-\tfrac{i}{v^*-\tau})^{-1} \\ \nonumber
& = & 1+\tfrac{i}{\ell_b}+O((\log T)^{-\frac{3}{2}}), \\ \label{expand3}
\tfrac{\tau-u^*}{\tau+i-u^*} &=  & 1-\tfrac{i_b}{\ell_b}+O((\log T)^{-\frac{3}{2}}).
\end{eqnarray}
By the law of the iterated logarithm, 
or by tail prob bounds of the standard normal,
\begin{eqnarray} \label{expand4}
|\wht \nu - \mu(\tau+1)| & = & O_p((\log T)^{-\frac{1}{3}}), \\ \label{expand5}
|\wht \mu - \mu(\tau)| & = & O_p((\log T)^{-\frac{1}{3}}).
\end{eqnarray}

By (\ref{711a})--(\ref{expand5}),
\begin{eqnarray*}
& & \sqrt{2 \ell_b}[Z_I(\tau+i)-Z_I(\tau)] \cr
& = & \lambda(\tau+i)[(\tfrac{v^*-\tau}{v^*-\tau-i}) \wht \nu - (\tfrac{\tau-u^*}{\tau+i-u^*}) \wht \mu
- \tfrac{2 \ell_b}{\lambda^2(\tau+i)}(S(\tau+i)-S(\tau))]-\lambda(\tau)(\wht \nu-\wht \mu) \cr
& = & {\rm (A)} + {\rm (B)} + {\rm (C)},
\end{eqnarray*}
where
\begin{eqnarray*}
{\rm (A)} & = & -\tfrac{2 \ell_b}{\lambda(\tau+i)} (S(\tau+i)-S(\tau)) \cr
& = & [1+O((\log T)^{-\frac{3}{2}})](2S^+(i)-2 i \Delta) \cr
& = & 2 S^+(i)-2i \Delta+(\log T)^{-\frac{3}{2}} O_p \Big( \max_{1 \leq i \leq (\log T)^{\frac{1}{4}}}
|S^+(i)-i \Delta| \Big) \cr
& = & 2 S^+(i) - 2i \Delta + O_p((\log T)^{-\frac{5}{4}}), \cr
{\rm (B)} & = & [\lambda(\tau+i)(\tfrac{v^*-\tau}{v^*-\tau-i})-\lambda(\tau)] \wht \nu \cr
& = & \tfrac{i \lambda(\tau+i)}{\ell_b} \wht \nu + O((\log T)^{-\frac{1}{2}}) \cr
& = & i \mu(\tau+1)+O_p((\log T)^{-\frac{1}{12}}), \cr
{\rm (C)} & = & -[\lambda(\tau+i)(\tfrac{\tau-\mu^*}{\tau+i-\mu^*})-\lambda(\tau)] \wht \mu \cr
& = & -i \mu(\tau)+O_p((\log T)^{-\frac{1}{12}}).
\end{eqnarray*}

We conclude (\ref{p+}) from the above expansions and $-2i \Delta+i(\mu(\tau+1)-\mu(\tau)) = -i \Delta$.
The calculations for (\ref{p-}),
with $-(\log T)^{\frac{1}{4}} \leq i \leq -1$, 
are similar.
$\wbox$

\subsubsection{Proof of (\ref{gscan})}

Select $\tau=\tau_j$ randomly from $\btau$ and let $\Delta=\Delta(\tau)$.
Let $A_T$ be defined as in (\ref{AT}) and let $\cI_T$ be all intervals $I$ satisfying (\ref{t/4}) and (\ref{show0}).  
By Lemma \ref{lem6} in Appendix C,
under $A_T$, 
with prob $\rightarrow 1$ the interval $I$ in the scan-CUSUM algorithm, 
for which $\tau \in I$, 
is in $\cI_T$.

Let $\wht \tau$ be the maximizer of $Z_I(t)$.
By (\ref{711a})--(\ref{expand5}),
$i^*=\wht \tau-\tau$ is asymptotically the mode of $\bp=(p(i): i \in \bZ)$ with
\begin{eqnarray} \label{pi2}
p(i) & = & \exp(\Delta S^+(i)-\tfrac{i \Delta^2}{2}), \\ \nonumber
p(-i) & = & \exp(\Delta S^-(i)-\tfrac{i \Delta^2}{2}), \\ \nonumber
p(0) & = & 1.
\end{eqnarray}
We conclude $\gamma_T(\Delta_0) \rightarrow \gamma_{\rm scan}(\Delta_0)$ by letting $\delta \rightarrow 0$.

\subsection{Proof of Theorem \ref{thm3}}

Let $\tau=\tau_j$ be a randomly selected change-point and let $k$ be such that $\tau \in I_k:=\{ 2k-1, 2k \}$.
Let $K = \lfloor \tfrac{T-1}{U} \rfloor$.
Let $J_k$ be the number of change-points in $I_k$.
Let $\Delta=\Delta(\tau)$.

Consider an oracle estimator
\begin{equation} \label{tau3}
\wht \tau_j = \left\{ \begin{array}{ll} \tau & \mbox{ if } J_{k-1}+J_k+J_{k+1} >1 \mbox{ or } k >K, \cr
{\rm arg \ max}_{t \in I_k} a(t) L_k(t) & \mbox{ otherwise,} \end{array} \right.
\end{equation}
where $L_k(t)$ is the likelihood defined in (\ref{Lkt}) with $U=2$.
Since $q \rightarrow 0$ the first case in (\ref{tau3}) occurs with prob$\rightarrow 0$.
In the analysis below we condition on the second case.

Since 
$$(a(t)|Y(t)=1) \sim G_q^* \mbox{ where } dG_q^*(a) = \tfrac{a}{q} dG_q(a),
$$
by (\ref{Ltt}),
the prob that $\wht \tau_j \neq \tau$ is asymptotically 
\begin{eqnarray*}
& & q^{-1} E_q(a(t) {\bf 1}_{\{ a(t) < a(t+1) \exp(\Delta S^+(1)-\frac{\Delta^2}{2}) \}}) \cr
& \geq & q^{-1} P(\Delta S^+(1)-\tfrac{\Delta^2}{2} > 0) E_q(a(t) {\bf 1}_{\{ a(t) < a(t+1) \}}),
\end{eqnarray*}
where $E_q$ is expectation under which $a(t)$ and $a(t+1)$ are i.i.d. $G_q$.

It is not possible for a change-point estimator to have smaller $\beta_T$ than the oracle estimator,
hence as $E_q(a(t) {\bf 1}_{\{ a(t) < a(t+1) \}}) = \tfrac{1}{2} E_q[\min(a(t),a(t+1))]$,
$$\limsup_{T \rightarrow \infty} (1-\beta_T) \geq 1-\tfrac{1}{2} P(Z > \tfrac{|\Delta|}{2}) \Big\{ 
\limsup_{q \rightarrow 0} q^{-1} E_q[\min(a(t),a(t+1))] \Big\}, 
$$
with $Z \sim {\rm N}(0,1)$ and $\Delta \sim F$ independent. 

By (\ref{VG}),
\begin{equation} \label{aa}
\limsup_{q \rightarrow 0} q^{-1} E_q[\min(a(t),a(t+1))] = \limsup_{q \rightarrow 0} q^{-1} \int_0^1 [1-G_q(x)]^2 dx > 0,
\end{equation}
and we conclude $\limsup_{T \rightarrow \infty} (1-\beta_T) > 0$.
$\wbox$

\subsection{Proof of Theorem \ref{thm4}}

Theorem \ref{thm4} follows from Lemmas \ref{lem3}--\ref{lem5}.
Lemmas \ref{lem3} and \ref{lem4} are proved in Sections 7.4.1 and 7.4.2.
Lemma \ref{lem5} is proved in Appendix D.

Let $E_q$ be expectation under which $a_0, a^+(i), a^-(i) \sim G_q$ for $i>0$,
$Z^+(u), Z^-(u) \sim {\rm N}(0,1)$ for $u>0$ and $\Delta \sim F$ are all independent.
Let $S^+(i) = \sum_{u=1}^i Z^+(u)$ and $S^-(i) = \sum_{u=1}^i Z^-(u)$.
Let
\begin{eqnarray*}
p(i) & = & \exp(\Delta S^+(i)-\tfrac{i \Delta^2}{2}), \cr
p(-i) & = & \exp(\Delta S^-(i)-\tfrac{i \Delta^2}{2}), \cr
B^+ & = & \{ a_0 \leq 2 a^+(i) p(i) \mbox{ for some } i >0 \}, \cr
B^- & = & \{ a_0 \leq 2 a^-(i) p(-i) \mbox{ for some } i >0 \}.
\end{eqnarray*}

Select $\tau=\tau_j$ randomly from $\btau$ and let $\Delta=\Delta(\tau)$.
Let $A_T$ be defined as in (\ref{AT}) and let $\cI_T$ be all intervals satisfying (\ref{t/4}) and (\ref{show0}).

\begin{lem} \label{lem3}
If $q^{-1} \int_0^1 [1-G_q(x)]^2 dx \rightarrow 0$ then $E_q(a_0 {\bf 1}_{B^+ \cup B^-}) \rightarrow 0$.
\end{lem}

\begin{lem} \label{lem4}
\begin{eqnarray} \label{lem4.1}
\sup \Big| \tfrac{1}{2}[Z_I^2(\tau+i)-Z_I^2(\tau)]-(\Delta S^+(i)-\tfrac{i \Delta^2}{2}) \Big| {\bf 1}_{A_T}
& \stackrel{p}{\rightarrow} & 0, \\ \label{lem4.2}
\sup \Big| \tfrac{1}{2}[Z_I^2(\tau-i)-Z_I^2(\tau)]-(\Delta S^-(i)-\tfrac{i \Delta^2}{2}) \Big| {\bf 1}_{A_T}
& \stackrel{p}{\rightarrow} & 0,
\end{eqnarray}
where the supremum is over $I \in \cI_T$ and $1 \leq i \leq (\log T)^{\frac{1}{4}}$.
\end{lem}

\begin{lem} \label{lem5}
$$\sum_{I \in \cI_T} \sum_{u \in I: |u-\tau| \geq (\log T)^{\frac{1}{4}}}
P(a(u) e^{\frac{1}{2} Z^2_I(u)} \geq a(\tau) e^{\frac{1}{2} Z_I^2(\tau)}, A_T) \rightarrow 0.
$$
\end{lem}

{\sc Proof of Theorem \ref{thm4}}.
Since $\wht \tau=\tau$ when $a(\tau+i) e^{\frac{1}{2} Z_I^2(\tau+i)} < a(\tau) e^{\frac{1}{2} Z_I^2(\tau)}$
for all $\tau+i \in I$ with $i \neq 0$,
by (\ref{PAT}), Lemmas \ref{lem4}--\ref{lem6} and the expansion of $p(i)$ in (\ref{pi2}),
\begin{equation} \label{1-beta}
\limsup (1-\beta_T) \leq \limsup P(a(\tau) \leq 2 a(\tau+i) p(i) \mbox{ for some } 1 \leq |i| \leq (\log T)^{\frac{1}{4}}).
\end{equation}
Let $a^+(i)=a(\tau+i)$ and $a^-(i)=a(\tau-1)$. 
Since
$$a(\tau) \sim G_q^* \mbox{ with } dG_q^*(a)=\tfrac{a}{q} dG_q(a),
$$
it follows from (\ref{1-beta}) that
$$\limsup (1-\beta_T) \leq q^{-1} E_q(a_0 {\bf 1}_{\{ B^+ \cup B^- \}}),
$$
and $\beta_T \rightarrow 1$ follows from Lemma \ref{lem3},
with $a_0[=a(\tau)] \sim G_q$ under $E_q$. $\wbox$

\subsubsection{Proof of Lemma \ref{lem3}}

Since $E_q(a_0 {\bf 1}_{B^+}) = E_q(a_0 {\bf 1}_{B^-})$,
it suffices to show 
\begin{equation} \label{pf3.1}
q^{-1} E_q(a_0 {\bf 1}_{B^+}|\Delta) \rightarrow 0 \mbox{ for all } \Delta \neq 0.
\end{equation}

Let $B^+(i) = \{ a_0 \leq 2a^+(i) p(i) \}$ and check that
\begin{eqnarray*} 
& & q^{-1} E_q(a_0 {\bf 1}_{B^+(i)}|\Delta) \cr
& = & q^{-1} E_q(a_0 {\bf 1}_{\{ a_0 \leq 2a^+(i) p(i), \Delta S^+(i) \leq \frac{i \Delta^2}{4} \}}|\Delta) \cr
& & \qquad + q^{-1} E_q(a_0 {\bf 1}_{\{ a_0 \leq 2a^+(i) p(i), \Delta S^+(i) > \frac{i \Delta^2}{4} \}}|\Delta) \cr
& \leq & q^{-1} E_q(2a^+(i) e^{-\frac{i \Delta^2}{4}}|\Delta) + q^{-1} E_q(a_0 {\bf 1}_{\{ \Delta S^+(i)
> \frac{i \Delta^2}{4} \}}|\Delta) \cr
& = & 2e^{-\frac{\Delta^2}{4}} + \Phi(-\tfrac{\sqrt{i} |\Delta|}{4}).
\end{eqnarray*}

Another useful bound is
\begin{eqnarray*} 
q^{-1} E_q(a_0 {\bf 1}_{B^+(i)}|\Delta) & \leq & q^{-1} E_q[\min(a_0, 2a^+(i) p(i))|\Delta] \\ \nonumber
& \leq & q^{-1} E_q[\min(a_0,a^+(i))(1+2 p(i))|\Delta] \\ \nonumber
& = & 3q^{-1} E_q[\min(a_0,a^+(i))].
\end{eqnarray*}

For any $i_0 > 0$,
\begin{eqnarray*}
& & q^{-1} E_q(a_0 {\bf 1}_{B^+}|\Delta) \cr
& \leq & \sum_{i=1}^{i_0} q^{-1} E_q(a_0 {\bf 1}_{B^+(i)}|\Delta) + \sum_{i=i_0+1}^{\infty} q^{-1} E_q(a_0 {\bf 1}_{B^+(i)}|\Delta)
\cr
& \leq & 3i_0 q^{-1} E_q[\min(a(1),a(2))] + \sum_{i=i_0+1}^{\infty} [2e^{-\frac{i \Delta^2}{4}}
+\Phi(-\tfrac{\sqrt{i} \Delta}{4})].
\end{eqnarray*}
Since $q^{-1} E_q[\min(a(1),a(2))] = q^{-1} \int_0^1 [1-G_q(x)]^2 dx \rightarrow 0$,
it follows that $\limsup_{q \rightarrow 0} q^{-1} E_q(a_0 {\bf 1}_{B^+}|\Delta)$ 
can be made arbitrarily small by choosing $i_0$ large enough.

\subsubsection{Proof of Lemma \ref{lem4}}

By (\ref{p+}) and (\ref{p-}) it suffices to show that
\begin{equation} \label{4/3}
\sup \Big| Z_I(\tau+i)+Z_I(\tau)-\Delta \sqrt{2 \ell_b} \Big| 
\Big| Z_I(\tau+i)-Z_I(\tau) \Big| {\bf 1}_{A_T} \stackrel{p}{\rightarrow} 0. 
\end{equation}

Recall $I=(t^*-\ell_b,t^*+\ell_b) = (u^*,v^*)$,
$\lambda(t) =\sqrt{(v^*-t)(t-u^*)}$ and $|\tau-t^*| \leq (\log T)^{\frac{1}{4}}$ for $I \in \cI_T$.
Since $v^*-u^* = 2 \ell_b$, 
by (\ref{ZI}),
$Z_I(\tau+i) \sim {\rm N}(\nu_i,1)$,
where 
\begin{eqnarray} \label{nui}
\nu_i & = & \left\{ \begin{array}{ll} \tfrac{\lambda(\tau+i)}{\sqrt{2 \ell_b}}
(1-\tfrac{i}{\tau+i-u^*}) \Delta & \mbox{ if } i > 0, \cr
\tfrac{\lambda(\tau+i)}{\sqrt{2 \ell_b}}
(1-\tfrac{|i|}{v^*-\tau+|i|}) \Delta & \mbox{ if } i < 0, \end{array} \right. \\ \nonumber
& = & [1+O((\log T)^{-\frac{3}{4}})] \sqrt{\tfrac{\ell_b}{2}} \Delta.
\end{eqnarray}
Since there are $O(\log T)$ intervals in $\cI_T$, 
it follows that $\sup |Z_I(\tau+i)-\nu_i| = O_p(\log \log T)$ and hence by (\ref{nui}),
\begin{equation} \label{Bon1}
\sup | Z_I(\tau)+Z_I(\tau+i)-\Delta \sqrt{2 \ell_b}| {\bf 1}_{A_T}
= O_p((\log T)^{-\frac{1}{4}}).
\end{equation}
By (\ref{p+}) and (\ref{p-}),
\begin{equation} \label{Bon2}
\sup |Z_I(\tau+i)-Z_I(\tau)| {\bf 1}_{A_T} = O_p((\log T)^{-\frac{1}{4}}).
\end{equation}
We conclude (\ref{4/3}) from (\ref{Bon1}) and (\ref{Bon2}).

\begin{appendix}

\section{Derivation of $\ell(a)$ and the EM procedure}

Given $\bX^n$, $1 \leq n \leq N$,
define latent variables $\btau = (\tau_j^n: 1 \leq n \leq N, 1 \leq j \leq \wht J^n)$,
with
$$P(\tau_j^n=t|\bX^n) = \tfrac{a(t) L_j^n(t)}{\sum_{u \in I^n(j)} a(u) L_j^n(u)} \mbox{ for } t \in I_j^n.
$$
The incomplete-data likelihood function,
based on a Poisson($\sum_{t=1}^{T-1} a(t)$) approximation for $J^n$,
is
$$L(\ba) = e^{\ell(a)} = e^{-N \sum_{t=1}^{T-1} a(t)} \prod_{n=1}^N \prod_{j=1}^{\hat J^n}
\Big( \sum_{t \in I^n(j)} a(t) L_j^n(t) \Big).
$$
The complete-data likelihood function is
$$L(\ba; \btau) = e^{-N \sum_{t=1}^{T-1} a(t)} \prod_{n=1}^N \prod_{j=1}^{\hat J^n} [a(t) L_j^n(t)]^{{\bf 1}_{\{ \tau_j^n=t \} }}.
$$

Given our current estimate of the parameters $\ba_k = (a_k(t): 1 \leq t \leq T-1)$,
the conditional prob of $\tau_j^n=t$ is
$$T_{jk}^n(t) = \tfrac{a_k(t) L_j^n(t)}{\sum_{u \in I^n_j} a_k(u) L_j^n(u)} \mbox{ for } t \in I_j^n.
$$
Hence the E-step,
\begin{eqnarray*}
Q(\ba|\ba_k) & = & E_{\btau|\ba_k}[\log L(\ba; \btau)] \cr
& = & -N \sum_{t=1}^{T-1} a(t) + \sum_{n=1}^N \sum_{j=1}^{\hat J^n} \sum_{t \in I^n(j)} T_{jk}^n(t)
\log(a(t) L_j^n(t)).
\end{eqnarray*}
The maximization of $Q(\ba|\ba_k)$ occurs when
$$- N + \sum_{(n,j): t \in I^n(j)} \tfrac{T_{jk}^n(t)}{a(t)} =0 
\Rightarrow a_{k+1}(t) = \tfrac{1}{N} \sum_{(n,j): t \in I^n(j)} T_{jk}^n(t).
$$

\section{Proofs of Lemmas \ref{lem1} and \ref{lem2}}

{\sc Proof of Lemma \ref{lem1}}. (a)
Let $s= \sum_{i=-\infty}^{\infty} r(i)$.
\begin{eqnarray*}
|{\rm med}(\br)| &= & \tfrac{1}{s} \sum_{i=-\infty}^{\infty} |{\rm med}(\br)| r(i) \cr
& \leq & \tfrac{1}{s} \sum_{i=-\infty}^{\infty} |i-{\rm med}(\br)| r(i) + \tfrac{1}{s} \sum_{i=-\infty}^{\infty} |i| r(i).
\end{eqnarray*}
Since $M(\br)={\rm med}(\br)$ minimizes $\sum_{i=-\infty}^{\infty} |i-M(\br)| r(i)$ over all $M(\br)$,
including $M(\br)=0$,
$$|{\rm med}(\br)| \leq \tfrac{2}{s} \sum_{i=-\infty}^{\infty} |i| r(i).
$$
Lemma \ref{lem1}(a) follows from $s \geq \max(r(0),r(i))$. 

\medskip
(b) Since $p(0)=1$,
\begin{equation} \label{Wp}
W(\bp) = \sum_{i=-\infty}^{\infty} |i| \min(1,p(i)),
\end{equation}
and so
$$E(W(\bp) {\bf 1}_{\{ |\Delta| \geq \Delta_0 \}}) = 2 \sum_{i=1}^{\infty}
i E[\min(1,e^{\Delta S^+(i)-\frac{i \Delta^2}{2}}) {\bf 1}_{\{ |\Delta| \geq \Delta_0 \}}].
$$
Since
\begin{eqnarray*}
E[\min(1,e^{\Delta S^+(i)-\frac{i \Delta^2}{2}})|\Delta] & = & \Phi(-\tfrac{\sqrt{i} |\Delta|}{2})
+ \int_{-\infty}^{\frac{i |\Delta|}{2}} \tfrac{1}{\sqrt{2 \pi i}} e^{-\frac{z^2}{2i}+\Delta z-\frac{i \Delta^2}{2}} dz \cr
& = & 2 \Phi(-\tfrac{\sqrt{i}|\Delta|}{2}),
\end{eqnarray*}
we conclude
$$E ( W(\bp) {\bf 1}_{\{ |\Delta| \geq \Delta_0 \}})
\leq 4 \sum_{i=1}^{\infty} i \Phi(-\tfrac{\sqrt{i} \Delta_0}{2}) < \infty. \mbox{ $\wbox$}
$$

\medskip
{\sc Proof of Lemma} \ref{lem2}. 
Let $s_T = \sum_{i=-\infty}^{\infty} r_T(i)$,  $s=\sum_{i=-\infty}^{\infty} r(i)$ and
$\delta_T = \sum_{i=-\infty}^{\infty} |r_T(i)-r(i)|$. 
By the triangle inequality, 
$|s_T-s| \leq \delta_T$.
Let
$$\epsilon = \min \Big( \sum_{t \leq {\rm med}(\br)} r(t)-\tfrac{s}{2},
\tfrac{s}{2}-\sum_{t < {\rm med}(\br)} r(t) \Big),
$$
which is positive by (\ref{lem2.1}).
For $T$ large such that $\delta_T \leq \frac{\epsilon}{2}$,
$$\min \Big( \sum_{t \leq {\rm med}(\br)} r_T(t)-\tfrac{s_T}{2},
\tfrac{s_T}{2}-\sum_{t < {\rm med}(\br)} r_T(t) \Big) \geq \epsilon - \tfrac{3}{2} \delta_T > 0,
$$
hence ${\rm med}(\br_T) = {\rm med}(\br)$.
$\wbox$

\section{Probability bounds for the proof of Theorem \ref{thm2}}

The probability bounds in Lemmas \ref{lem6} and \ref{lem7} are used for showing $\beta_T \rightarrow \beta_{\rm upper}$
and $\gamma_T(\Delta_0) \rightarrow \gamma_{\rm scan}(\Delta_0)$.
Recall
$$A_T = \{ \delta \leq |\Delta| \leq \delta^{-1}, \ Y(t)=0 \mbox{ for } 1 \leq |t-\tau| \leq 12 \delta^{-2} \log T \},
$$
where $\Delta=\Delta(\tau)$.

\begin{lem} \label{lem6}
{\rm (a)} For $6 \delta^{-2} \log T \leq \ell_b < 12 \delta^{-2} \log T$,
$$P(|Z_b(\tau)| < c_{\rm scan}, A_T) = o(\tfrac{1}{\log T}).
$$

\noindent {\rm (b)} For $18 \delta^{-2} \log T \leq \ell_b < 36 \delta^{-2} \log T$,
$$P(|Z_b(\tau)| < c_{\rm scan}, A_T) = o(\tfrac{1}{T}).
$$

\noindent {\rm (c)} 
$$\sum_{b: \ell_b \leq 2 \delta^2 \log T} \sum_{t: |t-\tau| < \ell_b} P(|Z_b(t)| \geq c_{\rm scan}, A_T) = o(\tfrac{1}{\log T}).
$$

\noindent {\rm (d)}
$$\sum_{b: 2 \delta^2 \log T \leq \ell_b \leq 12 \delta^{-2} \log T} \sum_{t: (\log T)^{\frac{1}{4}} < |t-\tau| < \ell_b}
P(|Z_b(t)| \geq |Z_b(\tau)|, A_T) = o(\tfrac{1}{\log T}).
$$
\end{lem}

{\sc Proof of Lemma} \ref{lem6}.
(a) Under $A_T$,
for $\ell_b \geq 6 \delta^{-2} \log T$,
$Z_b(\tau) \sim {\rm N}(\nu_b,1)$ with
$$\nu_b = \sqrt{\tfrac{\ell_b}{2}} |\Delta| \geq \sqrt{3 \log T}, 
$$
hence $P(|Z_b(\tau)| < c_{\rm scan}, A_T) \leq \Phi(-\sqrt{3 \log T}+c_{\rm scan}) = o(\tfrac{1}{\log T})$.

\medskip
(b) Under $A_T$,
for $\ell_b \geq 18 \delta^{-2} \log T$,
$Z_b(\tau) \sim {\rm N}(\nu_b,1)$ with
$$\nu_b = \sqrt{\tfrac{\ell_b}{2}} |\Delta| \geq 3 \sqrt{\log T}, 
$$
hence $P(|Z_b(\tau)| < c_{\rm scan}, A_T) \leq \Phi(-3 \sqrt{\log T}+c_{\rm scan}) = o(\tfrac{1}{T})$.

\medskip
(c) Under $A_T$,
for $\ell_b \leq 2 \delta^2 \log T$ and $\tau-\ell_b < t < \tau+\ell_b$,
$Z_b(t)$ is unit variance normal with
$$|E Z_b(t)| \leq |E Z_b(\tau)| \leq \sqrt{\log T}, 
$$
hence the sum of prob in Lemma \ref{lem6}(c) is $O((\log T)^2) \Phi(-c_{\rm scan}+\sqrt{\log T}) = o(\tfrac{1}{\log T})$.

\medskip
(d) Under $A_T$,
$Z_b(t)-Z_b(\tau) \sim {\rm N}(-\tfrac{|t-\tau|}{\sqrt{2 \ell_b}} \Delta, \tfrac{|t-\tau|}{\ell_b})$ and
$Z_b(\tau) \sim {\rm N}(\Delta \sqrt{\tfrac{\ell_b}{2}},1)$,
hence the sum of prob in Lemma \ref{lem6}(d) is
$$O(\log T) \Big[ \sum_{i > (\log T)^{\frac{1}{4}}} \Phi \Big( -\delta \sqrt{\tfrac{i}{2}} \Big)
+ \Phi(-\delta^2 \log T) \Big]
=o(\tfrac{1}{\log T}). \mbox{ $\wbox$}
$$

\medskip
Recall 
$$\cI_T = \{ (t^*-\ell_b,t^*+\ell_b): |t^*-\tau| \leq (\log T)^{\frac{1}{4}}, 2 \delta^2 \log T < \ell_b < 12 \delta^{-2}
\log T \}.
$$

\begin{lem} \label{lem7}
$$\sum_{I \in \cI_T} \sum_{t \in I: |t-\tau| \geq (\log T)^{\frac{1}{4}}}
|t-\tau| P(|Z_I(t)| \geq |Z_I(\tau)|, A_T) \rightarrow 0.
$$
\end{lem}

{\sc Proof}.
Let $I=(u^*,v^*)$ and recall $\lambda(t)=\sqrt{(v^*-t)(t-u^*)} = \sqrt{\ell_b^2-|t^*-\tau|^2}$.
Consider $\Delta > 0$.
Under $A_T$, for $t \in I$ with $t>\tau$,
\begin{eqnarray} \label{EZ2}
E(Z_I(t)-Z_I(\tau)) & = & [\tfrac{\lambda(t)}{\sqrt{2 \ell_b}}(1-\tfrac{t-\tau}{t-u^*})-\tfrac{\lambda(\tau)}{\sqrt{2 \ell_b}}]
\Delta \\ \nonumber
& = & -\tfrac{\lambda(\tau) \Delta}{\sqrt{2 \ell_b}} \Big[ 1 -\sqrt{\tfrac{(v^*-t)(\tau-u^*)}{(v^*-\tau)(t-u^*)}} \Big].
\end{eqnarray}
Moreover
\begin{eqnarray*}
{\rm Cov}(Z_I(t),Z_I(\tau)) & = & \tfrac{\lambda(t) \lambda(\tau)}{2 \ell_b}
(\tfrac{1}{t-u^*}+\tfrac{1}{v^*-\tau}-\tfrac{t-\tau}{(t-u^*)(v^*-\tau)}) \cr
& = & \sqrt{\tfrac{(v^*-t)(\tau-u^*)}{(v^*-\tau)(t-u^*)}},
\end{eqnarray*}
hence
\begin{equation} \label{v7}
{\rm Var}(Z_I(t)-Z_I(\tau)) = 2 \Big[1-\sqrt{\tfrac{(v^*-t)(\tau-u^*)}{(v^*-\tau)(t-u^*)}} \Big].
\end{equation}

For $|t-\tau| \geq (\log T)^{\frac{1}{4}}$ and $|t^*-\tau| \leq (\log T)^{\frac{1}{4}}$,
\begin{eqnarray} \label{v7.1}
\sqrt{\tfrac{(v^*-t)(\tau-u^*)}{(v^*-\tau)(t-u^*)}} & = & (1-\tfrac{t-\tau}{v^*-\tau})^{\frac{1}{2}}
(1-\tfrac{t-\tau}{t-u^*})^{\frac{1}{2}} \\ \nonumber
& \leq & (1-[1+o(1)] \tfrac{(\log T)^{\frac{1}{4}}}{2 \ell_b})^2 \\ \nonumber
& = & 1-[1+o(1)] \tfrac{(\log T)^{\frac{1}{4}}}{\ell_b}.
\end{eqnarray}

By (\ref{EZ2})--(\ref{v7.1}),
\begin{eqnarray} \label{ZItt}
\tfrac{E[Z_I(t)-Z_I(\tau)]}{\sqrt{{\rm Var}(Z_I(t)-Z_I(\tau))}}
& = & -\tfrac{\lambda(\tau) \Delta}{2 \sqrt{\ell_b}}
\Big[ 1 -\sqrt{\tfrac{(v^*-t)(\tau-u^*)}{(v^*-\tau)(t-u^*)}} \Big]^{\frac{1}{2}} \\ \nonumber
& \leq & -[\Delta+o(1)] \tfrac{\lambda(\tau) (\log T)^{\frac{1}{8}}}{2 \ell_b} \\ \nonumber
& \leq & -[\tfrac{\Delta}{2}+o(1)] (\log T)^{\frac{1}{8}}.
\end{eqnarray}
The calculations are similar for $\Delta < 0$.

Hence for $T$ large,
$$P(|Z_I(t)| \geq |Z_I(\tau)|, A_T) \leq \Phi(-\tfrac{\delta (\log T)^{\frac{1}{8}}}{3})
=o(\tfrac{1}{(\log T)^3}),
$$
and Lemma \ref{lem7} follows from $\# \cI_T = O(\log T)$.
$\wbox$

\section{Proof of Lemma \ref{lem5}}

Let $u \in I$ with $I \in \cI_T$ and $|u-\tau| \geq (\log T)^{\frac{1}{4}}$.
Let 
$$W_I(u,\tau) = \tfrac{Z_I(u)-Z_I(\tau)-E[Z_I(u)-Z_I(\tau)]}{\sqrt{{\rm Var}(Z_I(u)-Z_I(\tau))}} (\sim {\rm N}(0,1)).
$$
Express
\begin{eqnarray} \label{8.1}
& & \tfrac{1}{2}(Z_I^2(u)-Z_I^2(\tau)) \\ \nonumber
& = & \{ \tfrac{E[Z_I(u)-Z_I(t)]}{\sqrt{{\rm Var}(Z_I(u)-Z_I(t))}}  + W_I(u,\tau) \} \\ \nonumber
& & \qquad \times \tfrac{1}{2}[Z_I(u)+Z_I(\tau)] \sqrt{{\rm Var}(Z_I(u)-Z_I(\tau))}.
\end{eqnarray}

Consider $\Delta>0$.
Since $I \in \cI_T$,
$|\tau-t^*| \leq (\log T)^{\frac{1}{4}}$ and $\ell_b > 2 \delta^2 \log T$,
see (\ref{t/4}) and (\ref{show0}),
\begin{equation} \label{lambdatau}
\lambda(\tau) \sim \sqrt{\ell_b^2-|\tau-t^*|^2} \sim \ell_b.
\end{equation}
Let $W_I(u) = Z_I(u)-E(Z_I(u)) (\sim {\rm N}(0,1))$.
By (\ref{lambdatau}) and $E[Z_I(u)] > 0$, 
for $|\tau-t^*| \leq (\log T)^{\frac{1}{4}}$,
\begin{eqnarray} \label{8.2}
\tfrac{1}{2}[Z_I(u)+Z_I(\tau)] & > & \tfrac{1}{2} E(Z_I(\tau))+W_I(u)+W_I(\tau) \\ \nonumber
& = & \tfrac{\lambda(\tau) \Delta}{2 \sqrt{2 \ell_b}} + W_I(u)+W_I(\tau) \\ \nonumber
& = & [\Delta+o(1)] \sqrt{\tfrac{\ell_b}{8}} + W_I(u)+W_I(\tau).
\end{eqnarray}

By (\ref{v7}) and (\ref{v7.1}),
\begin{equation} \label{8.3}
\sqrt{{\rm Var}(Z_I(u)-Z_I(\tau))} \sim (\log T)^{\frac{1}{8}} \sqrt{\tfrac{2}{\ell_b}}.
\end{equation}

Substituting (\ref{ZItt}), (\ref{8.2}) and (\ref{8.3}) into (\ref{8.1}) gives us
\begin{eqnarray*}
\tfrac{1}{2}(Z_I^2(u)-Z_I^2(\tau)) & \leq & -[\tfrac{\Delta^2}{4}+o(1)] (\log T)^{\frac{1}{4}} \cr
& & \qquad + [W_I(u,\tau)-W_I(u)-W_I(\tau)][\tfrac{\Delta}{2}+o(1)] (\log T)^{\frac{1}{8}}.
\end{eqnarray*}
The calculations are similar for $\Delta<0$.
Hence
$$P(\tfrac{1}{2} [Z_I^2(u)-Z_I^2(\tau)] \geq -\tfrac{\delta^2}{5} (\log T)^{\frac{1}{4}}, 
\ |\Delta| \geq \delta) = o(\tfrac{1}{(\log T)^3}).
$$
Hence it suffices to show that 
\begin{equation} \label{8.4}
P(a(u) \geq a(\tau) \exp(\tfrac{\delta^2}{5} (\log T)^{\frac{1}{4}}), \ |\Delta| \geq \delta) = o(\tfrac{1}{(\log T)^3}).
\end{equation}

Indeed as $a(\tau) \sim G_q^*$ (with $dG_q^*(a) = \tfrac{a}{q} dG_q(a)$) and $a(u) \sim G_q$,
the prob in (\ref{8.4}) is bounded by
\begin{eqnarray*}
& & q^{-1} E_q(a(\tau) {\bf 1}_{\{ a(u) \geq a(\tau) \exp(\frac{\delta^2}{5} (\log T)^{\frac{1}{4}}) \}}) \cr
& \leq & q^{-1} \exp(-\tfrac{\delta^2}{5} (\log T)^{\frac{1}{4}}) E_q(a(u)) \cr
& = & \exp(-\tfrac{\delta^2}{5} (\log T)^{\frac{1}{4}}) = o(\tfrac{1}{\log T}). 
\end{eqnarray*}
\end{appendix}


\begin{thebibliography}{1}
\bibitem{ADH05} \textsc{Arias-Castro, E., Donoho, D.} and \textsc{Huo X.} (2005).
Near-optimal detection of geometric objects by fast multiscale methods.
\textit{IEEE Trans. Inf. Theory} \textbf{51} 2402--2425.

\bibitem{ADH06} \textsc{Arias-Castro, E., Donoho, D.} and \textsc{Huo X.} (2006). 
Adaptive multiscale detection of filamentary structures in a background of uniform noise.
\textit{Ann. Statist.} \textbf{34} 326--349. 

\bibitem{BH92} \textsc{Barry, D.} and \textsc{Hartigan, J.A.} (1993).
A Bayesian analysis for change point models.
\textit{J. Amer. Statist. Assoc.} \textbf{88} 309--319. 


\bibitem{CW15} \textsc{Chan, H.P.} and \textsc{Walther, G.} (2015). 
Optimal detection of multi-sample aligned sparse signals.
\textit{Ann. Statist.} \textbf{43} 1865--1895


\bibitem{Chi98} \textsc{Chib, S.} (1998).
Estimation and comparison of multiple change-point models.
\textit{J. Econometrics} \textbf{86} 221--241.

\bibitem{Cho16} \textsc{Cho, H.} (2016).
Change-point detection in panel data via double CUSUM statistic. 
\textit{Electron. J. Statist.} {\bf 10} 2000--2038.

\bibitem{CF15} \textsc{Cho, H.} and \textsc{Fryzlewicz, P.} (2015).
Multiple change-point detection for high-dimensional time series via sparsified binary segmentation. \textit{J. Roy. Statist. Soc.} B \textbf{77} 475--507.


\bibitem{DKK16}
\textsc{Du, C., Kao, C.L.} and \textsc{Kou, S.C.} (2016).
Stepwise signal extraction via marginal likelihood.
\textit{J. Amer. Statist. Assoc.} {\bf 111} 314--330.


\bibitem{EH19} \textsc{Enikeeva, F.} and \textsc{Harchaoui, Z.} (2019).
High-dimensional change-point detection under sparse alternatives.
\textit{Ann. Statist.} {\bf 47} 2051--2079.


\bibitem{Fry14} \textsc{Fryzlewicz, P.} (2014).
Wild binary segmentation for multiple change-point detection.
\textit{Ann. Statist.} {\bf 42} 2243--2281.


\bibitem{HH14} \textsc{Horv\'{a}th L.} and \textsc{Hu\v{s}kov\'{a} M.} (2014).
Change-point detection in panel data.
\textit{J. Time Ser. Anal.} {\bf 23} 631--648.



\bibitem{JCL13} \textsc{Jeng, X. J.}, \textsc{Cai, T. T.} and \textsc{Li, H.} (2013).
Simultaneous discovery of rare and common segment variants. 
\textit{Biometrika} {\bf 100} 157--172.

\bibitem{Jir15} \textsc{Jirak, M.} (2015).
Uniform change point tests in high dimension.
\textit{Ann. Statist.} {\bf 43} 2451--2483.

\bibitem{KCG15} \textsc{Ko, S.}, \textsc{Chong, T.} and \textsc{Ghosh, P.} (2015).
Dirichlet process hidden markov multiple change-point model.
\textit{Bayesian Anal.} {\bf 10} 275--296.

\bibitem{LX11} \textsc{Lai, T.L.} and \textsc{Xing, H.} (2011).
A simple Bayesian approach to multiple change-points.
\textit{Statist. Sinica} {\bf 21} 539--569.

\bibitem{LGS21} \textsc{Liu, H.}, \textsc{Gao, C.} and \textsc{Samworth, R.} (2021).
Minimax rate in sparse high-dimensional change-point detection.
\textit{Ann. Statist.} {\bf 49} 1081--1112.

\bibitem{Mei10} \textsc{Mei, Y.} (2010).
Efficient scalable schemes for monitoring large number of data streams.
{\it Biometrika} {\bf 97} 419--433.


\bibitem{NZ12} \textsc{Niu, Y.S.} and \textsc{Zhang, H.} (2012).
The screening and ranking algorithm to detect DNA copy number variation.
\textit{Ann. Appl. Statist.} {\bf 6} 1306--1326.


\bibitem{PCV20} \textsc{Pilliat, E., Carpentier. A.} and \textsc{Verzelen, N.} (2020).
Optimal multiple change-point detection for high-dimensional data.
arXiv preprint {\it arXiv:2011.07818}.

\bibitem{Sie85} \textsc{Siegmund, D.} (1985).
{\it Sequential Analysis: Tests and Confidence Intervals}. Springer, New York. 


\bibitem{WS18} \textsc{Wang, T.} and \textsc{Samworth, R.} (2018).
High dimensional change point estimation via sparse projection.
\textit{J. Roy. Statist. Soc.} B {\bf 80} 57--83.


\bibitem{XS13} \textsc{Xie, Y.} and \textsc{Siegmund, D.}
(2013). Sequential multi-sensor change-point detection.
\textit{Ann. Statist.} \textbf{41} 670--692.

\bibitem{Yao84} \textsc{Yao, Y.C.} (1984).
Estimation of a noisy discrete-time step function: Bayes and empirical Bayes approaches. \textit{Ann. Statist.} {\bf 12} 1434--1447.

\bibitem{ZS07} \textsc{Zhang, N.R.} and \textsc{Siegmund, D.} (2007).
A modified Bayes information criterion with applications to the analysis of comparative genomic hybridization.
\textit{Biometrics} \textbf{63} 22--52.

\bibitem{ZSJL10} \textsc{Zhang, N.R., Siegmund, D., Ji, H.} and \textsc{Li, J.Z.}
(2010). Detecting simultaneous changepoints in multiple sequences.
\textit{Biometrika} \textbf{97} 631--645.
\end{thebibliography}
\end{document}